\documentclass[12pt]{amsart}
\usepackage[active]{srcltx}%%%% Inverse search
\usepackage{amsmath,amstext,amsthm,amsfonts,latexsym}
\usepackage{graphicx}
\usepackage{graphics}
\usepackage[all,knot,poly]{xy}
\xyoption{arc}
\usepackage{a4wide}
\usepackage[latin1]{inputenc}
\usepackage[colorlinks=true,linkcolor=red,urlcolor=black]{hyperref}
\usepackage[pagewise]{lineno}

\numberwithin{equation}{section}

\newtheorem{teo}{Theorem}[section]
\newtheorem{prop}[teo]{Proposition}

\newtheorem{lema}[teo]{Lemma}

\newtheorem{maintheorem}{Theorem}
\newtheorem{maincor}[maintheorem]{Corollary}
\newtheorem{mainlemma}{Lemma}

\newcommand{\bmt}{\begin{maintheorem}}
\newcommand{\emt}{\end{maintheorem}}
\newcommand{\bml}{\begin{mainlemma}}
\newcommand{\eml}{\end{mainlemma}}
\numberwithin{equation}{section}

\newcommand{\La}{\Lambda}

\newcommand{\e}{\varepsilon}

\newcommand{\si}{\sigma}

\newcommand{\de}{\delta}

\newcommand{\B}{\mathcal{B}}

\newcommand{\supp}{{\rm supp}\:}

\renewcommand{\v}{{\bf v}}

\title{On the non-robustness of intermingled basins}
\date{\today}

\author[R. Ures]{Ra\'ul Ures}
\address{Ra\'ul Ures,IMERL, Facultad de Ingenier\'ia, Universidad de la Rep\'ublica, CC 30, Montevideo-Uruguay.}
\email{ures@fing.edu.uy}

\author[C. H. V\'asquez]{Carlos H. V\'asquez}
\address{Carlos H. V\'asquez, Instituto de Matem\'atica,
Pontificia Universidad Cat\'olica de Valpara\'{\i}so, Blanco Viel 596,
Cerro Bar\'on, Valpara\'{\i}so-Chile.} \email{carlos.vasquez@ucv.cl}
\begin{thanks}{CV  was supported by the Center of Dynamical Systems and Related
Fields c\'odigo ACT1103 PIA - Conicyt and Proyecto Fondecyt 1130547.}
\end{thanks}


\begin{document}

\begin{abstract}
It is well-known that it is possible to construct a partially hyperbolic diffeomorphism on the 3-torus in a similar way than in Kan's example. It has two hyperbolic physical measures with intermingled basins supported on two embedded tori with Anosov dynamics. A natural question is how robust is the intermingled basins phenomenon for diffeomorphisms defined on boundaryless manifolds? In this work we  study  partially hyperbolic diffeomorphisms on the 3-torus and show that the intermingled basins phenomenon is not robust.
\end{abstract}

\maketitle

%%%%%%%%%%%%%%%%%%%%%%%%%%%%%%%%%%%%%%%%%%%%%%%%%%%%%%%%%%%%%%%

\section{Introduction}\label{sec:intro}

Attractors play a key role in the study of non-conservative dynamics. The description of attractors and the properties  of their basins help predict the future behaviour of the orbits of a system. In this work we deal with physical measures i.e. an ergodic measure $\mu$ is  physical if its basin of attraction has positive volume (see Section \ref{sec:Preliminaries} for precise definitions). We will think these measures as the attractors of our systems. 

In many cases, basins are (essentially) open sets and it is clear that if a point belongs to certain regions its trajectory goes, almost surely,  to an attractor that is well determined. 
For instance, uniformly hyperbolic diffeomorphisms exhibit a finite number of physical measures and the union of their  basins cover Lebesgue almost every point the ambient manifold. Moreover, each one of their basins is an open set  (modulo a set of null volume) and then, we can clearly distinguish  one attractor from the others.

Outside  the uniformly hyperbolic world, this kind of behaviour of the basins of attractors is no longer true. Open sets of diffeomorphisms of manifolds with boundary may have attractors with intermingled basins.  More specifically, two or more basins are dense in the same open set. It was I. Kan~\cite{Kan:1994kw}  (See also \cite{BDV} for a description of the example in terms of the partial hyperbolicty and Lyapunov exponents) who showed for the first time the existence of  examples of  partially hyperbolic endomorphisms defined on a surface and exhibiting two hyperbolic physical measures whose basins are intermingled. Moreover,  he showed that such phenomenon is robust among the maps preserving the boundary. We refer the reader to \cite{IKS08} for a rigorous proof of Kan example and \cite{BM2008} for a generalization of the Kan example and its relation with the sign of the Schwarzian derivative. In \cite{KS2011} the authors  shown that the set of points that are not attracted by either of the components in the Kan's example has Hausdorff dimension less than the dimension of the phase space itself.  Following the same type of arguments, it is possible to construct a partially hyperbolic diffeomorphism defined on a 3-manifold with boundary exhibiting two intermingled physical measures, and such phenomenon still can be made robust. Furthermore, it is well known that it is possible to extend such example to the 3-torus, but  in this case it is no longer robust. We describe these examples in Section~\ref{sec:Examples}.

The existence of these examples rise the question of how robust are the intermingled basins phenomenon for diffeomorphisms defined on boundaryless manifolds. In this work we  show that   partially hyperbolic diffeomorphisms on the 3-torus having hyperbolic physical measures with intermingled basins are not robust. 

In a recent work, Okunev \cite{Okunev}, studied attractors in the sense of Milnor  in the most restrictive case of $C^r$ partially hyperbolic skew products on $\mathbb{T}^3$ with an Anosov dffeomorphisms acting on the base $\mathbb{T}^2$. The author obtains results with the same flavour as ours without any explicit hypotheses about Lyapunov exponent in the central direction.

%\section{Setting and Results}\label{sec:results}
We are interested in diffeomorphisms defined on a 3-dimensional manifold $M$, in particular we put our focus on $M=\mathbb{T}^3$. We give some basic definitions necessary to formulate the results, but the reader can find the precise definitions, properties and more detailed information in Section~\ref{sec:Preliminaries} and the references therein. 

A diffeomorphism $f:M\to M$ is \emph{partially hyperbolic} if the tangent bundle splits into three non trivial sub-bundles
$$TM =E^{uu}\oplus E^c \oplus E^{ss}$$
such that the strong stable sub-bundle $E^{ss}$ is uniformly contracted, the strong unstable sub-bundle $E^{uu}$ is uniformly expanded and
 the center sub-bundle $E^c$ may contract or expand, but this contractions or expansions are weaker than the strong expansions and contractions of the corresponding strong sub-bundles. 

%the sub-bundles $E^u$, $E^c$, and $E^s$ are  H\"{o}lder continuous and that 
It is known that 
there are unique foliations $W^{uu}$ and $W^{ss}$ tangent to $E^{uu}$ and $E^{ss}$ respectively \cite{BP1974,HPS} but in general, $E^c$, $E^{cu} = E^c \oplus E^{uu}$, and $E^{cs} = E^c \oplus E^{ss}$ do not integrate to foliations (see \cite{HHU}). 
The system is said to be \emph{dynamically coherent} if there exist invariant foliations $W^{cu}$ and $W^{cs}$ tangent to  $E^{cu}$ and  $E^{cs}$ respectively.  Of course, if this is the case, there exists an invariant foliation tangent to $E^c$ obtained just by intersecting $W^{cu}$ and $W^{cs}$. 
We will study dynamically coherent diffeomorphism with compact center leaves. As we mentioned above these diffeomorphisms are not always dynamically coherent although there are some results providing this property. Just to mention one result,  Brin, Burago, and Ivanov have shown that every absolute partially hyperbolic system (see Subsection  \ref{ssec:ph} for the definition) on the 3-torus is dynamically coherent \cite{BBI2009}.

A set $K\subseteq M$ is \emph{$u$-saturated} if it is the union of complete strong unstable leaves. The diffeomorphism $f$ is \emph{accessible} if 
every pair of points $x,y\in M$ can be joined by an arc consisting of finitely many segments contained in the leaves of the strong stable and strong unstable foliations. Assuming that the center bundle is one-dimensional, K. Burns, F. R. Hertz, J. R. Hertz, A. Talitskaya and R. Ures \cite{BHHTU} proved that the accessibility property is open and dense among the $C^r$-partially hyperbolic diffeomorphisms (see also \cite{NT2001}) . Our main theorem is the following.

\begin{maintheorem}\label{mTeo:A}
Let  $f\in{\rm Diff}^r(\mathbb{T}^3)$, $r\geq 2$, be partially hyperbolic, dynamically coherent with compact center leaves. Let $\mu$ be a physical measure with negative center Lyapunov exponent. Assume that $K\subseteq \mathbb{T}^3$ is  a compact,  $f$-invariant and $u$-saturated subset such that $K\subseteq \overline{\mathcal{B(\mu)}}\setminus\supp\mu$. Then, $K$ contains a finite union of periodic  2-dimensional $C^1$-tori, tangent to $E^u\oplus E^s$. In particular $f$ is not accessible.
\end{maintheorem}

We say that  two physical measures $\mu$ and $\nu$ with disjoint supports have {\em intermingled basins} \cite{Kan:1994kw} if  for an open set $U\subseteq M$ we have
${\rm Leb}(V\cap \mathcal{B}(\mu))>0$ and  ${\rm Leb}(V\cap \mathcal{B}(\nu))>0$ for any open set $V\subset U$.

\begin{maincor}\label{mcor:B}
The set of dynamically coherent partially hyperbolic  $C^r$-diffeomorphisms defined on $\mathbb{T}^3$, $r\geq 2$, exhibiting intermingled hyperbolic physical measures has empty interior. 

Moreover, if $f:\mathbb{T}^3\to\mathbb{T}^3$  is isotopic to a hyperbolic automorphism, there do not exist hyperbolic physical measures with intermingled basins.
\end{maincor}

Closely related, Hammerlindl and Potrie \cite{HP2014} showed that partially hyperbolic diffeomorphisms on $3$-nilmanifold admit a unique $u$-saturated minimal subset. Then, $f$ has a unique hyperbolic physical measure (see Section~\ref{ssec:metric} for more details) and thus,  it is not possible to have the intermingled basins phenomenon. We have as corollary of their work:

\begin{maincor}\label{mcor:C} If $M$ is a $3$-nilmanifold , then there does not exist hyperbolic physical measures with intermingled basins.
\end{maincor}

This paper is organized as follows. Section~\ref{sec:Preliminaries} is devoted to introduce the main tools in the proof: partial hyperbolic diffeomorphisms, physical measures, $u$-measures and Lyapunov exponents. A toy example as well as  Kan-like examples are revisited in Section~\ref{sec:Examples}. Proofs of Theorem~\ref{mTeo:A} and Corollary~\ref{mcor:B} are developed in Section~\ref{sec:proofs}.

\section{Preliminaries}\label{sec:Preliminaries}

\subsection{Partial hyperbolicity}\label{ssec:ph}

 Throughout this paper we shall work with a {\em
partially hyperbolic diffeomorphism} $f$, that is, a diffeomorphism
admitting a nontrivial $Tf$-invariant splitting of the tangent
bundle $TM = E^{ss}\oplus E^c \oplus E^{uu}$, such that all unit vectors
$v^\si\in E^\si_x$ ($\si= ss, c, uu$) with $x\in M$ satisfy:
$$\|T_xfv^{ss}\| < \|T_xfv^c\| < \|T_xfv^{uu}\| $$
for some suitable Riemannian metric. $f$ also must satisfy
that
$\|Tf|_{E^{ss}}\| < 1$ and $\|Tf^{-1}|_{E^{uu}}\| < 1$.  We also want to introduce a stronger type of partial hyperbolicity. We will say that $f$ is {\em absolutely partially hyperbolic}\, if it is partially hyperbolic and $$\|T_xfv^{ss}\| < \|T_yfv^c\| < \|T_zfv^{uu}\| $$ for all $x,y,z\in M$.\par%

 For partially hyperbolic diffeomorphisms, it is a well-known fact that there are foliations $W^\si$ tangent to the distributions $E^\si$
for $\si=ss,uu$ . The leaf of $W^{\si}$
containing $x$ will be called $W^{\si}(x)$, for $\si=ss,uu$. \par

In general it is not true that there is a foliation tangent to
$E^c$. Sometimes there is no foliation tangent to $E^c$. Indeed, there may be no foliation tangent to $E^c$ even if  $\dim E^c =1$ (see \cite{HHU}). We shall say that $f$ is {\em dynamically coherent} if there exist invariant foliations $W^{c\si}$ tangent to $E^{c\si}=E^c \oplus E^\si$  for $\si=ss,uu$. Note that by taking the intersection of these foliations we obtain an invariant foliation $W^c$ tangent to $E^c$ that subfoliates $W^{c\si}$ for $\si=s,u$. In this paper all partially hyperbolic diffeomorphisms will be dynamically coherent.

 We shall say that a set $X$ is {\em $\si$-saturated} if it is a union
of leaves of the strong foliations $\mathcal{W}^\si$ for $\si=ss$ or $uu$. We also say
that $X$ is $su$-saturated if it is both $s$- and $u$-saturated. The
 accessibility class of the point $x\in M$ is the
minimal $su$-saturated set containing $x$.  In case there is some
$x\in M$ whose accessibility class is $M$, then the diffeomorphism
$f$ is said to have the {\em accessibility property}. This is
equivalent to say that any two points of $M$ can be joined by a path
which is piecewise
tangent to $E^{ss}$ or to $E^{uu}$. \par%

\subsection{Physical measures, $u$-measures, Lyapunov exponents}\label{ssec:metric}

In this section we consider $f\colon M\rightarrow M$ be a diffeomorphism, not necessarily partially hyperbolic, defined on the riemannian manifold $M$. We denote by ${\rm Leb}$ the normalized volume form on $M$. 

A point $z\in M$ is {\em Birkhoff regular} if  the Birkhoff averages
\begin{equation}\label{eq:birkhoffneg}
 \varphi^-(z)=\lim_{n\to\infty}\frac{1}{n}\sum_{k=0}^{n-1} \varphi(f^{-k}(z)),
\end{equation}

\begin{equation}\label{eq:birkhoffpos}
 \varphi^+(z)=\lim_{n\to\infty}\frac{1}{n}\sum_{k=0}^{n-1} \varphi(f^k(z));
\end{equation}
are defined and $\varphi^-(z)=\varphi^+(z)$ for every $\varphi:M\to\mathbb{R}$ continuous. We denote by 
${\mathcal R}(f)$ the set of Birkhoff regular points of $f$. Birkhoff Ergodic Theorem \cite{M87, W82}, implies that the set ${\mathcal R}(f)$ has full measure with respect to any $f$-invariant 
measure $\xi$. When $\xi$ is an ergodic measure,
$$\varphi^-(z)=\varphi^+(z)=\int_M  \varphi\:d\xi,$$
for every $z$ in a $\xi$-full measure set ${\mathcal R}(\xi)$. 

If $\xi$ is an $f$-invariant measure, the {\em basin} of $\xi$ is the set

$${\mathcal B}(\xi)=\{z\in M\::\: \lim_{n\to\infty}\frac{1}{n}\sum_{k=0}^{n-1} \varphi(f^k(z))=\int_M  \varphi\:d\xi, \mbox{ for all } \varphi\in C^0(M,\mathbb{R})\}$$

If $\xi$ is an $f$-invariant ergodic measure, then  ${\mathcal R}(\xi)\subseteq{\mathcal B}(\xi)$,and so ${\mathcal B}(\xi)$ has full $\xi$-measure. 

An $f$-invariant probability measure $\mu$ is \emph{physical} if its basin
${\mathcal B}(\mu)$
has positive Lebesgue measure on $M$ \cite{BDV,Young:2002vc}. A physical measure is said to be \emph{hyperbolic} if all its Lyapunov exponents are nonzero \cite{BP}. In the setting of partially hyperbolic diffeomorphims defined on a 3-dimensional manifold, a physical measure is hyperbolic if 

$$\lambda^c(\mu)=\int\log\|Df|E^c\|d\mu\ne 0.$$

A point $x\in M$ is {\em Lyapunov regular} if there exist an integer $p(x)\leq \dim M$, numbers
$$\chi_1(x)<\dots<\chi_{p(x)}(x),$$
and a decomposition
\begin{equation}\label{eq:Osedecom}
T_xM=\bigoplus_{i=1}^{p(x)} H_i(x)
\end{equation}
into subspaces $H_i(x)$ such that $Df(x)H_j(x)=H_j(f(x))$, and for every $v\in H_j(x)\setminus\{0\}$
\begin{equation}\label{eq:charval}
\chi_j(x)=\lim_{n\to\pm \infty}\frac{1}{n}\log\|Df^n(x)v\|.
\end{equation}
Denote by $\Lambda(f)$ the set of Lyapunov regular points. The numbers $\chi_1(\xi)<\dots<\chi_{p(x)}(x),$ are called the {\em Lyapunov exponents} of $x$. The splitting  \eqref{eq:Osedecom} is called {\em Oseledets decomposition} and the subspaces $H_i(x)$ are called {\em Oseledets subespaces} at $x$. Oseledet's Theorem \cite{Oseledec:1968tk,M87} guarantee that the set $\Lambda(f)$ has full measure with respect any invariant measure. In general the functions $x\to\chi_j(x)$, $x\to H_j(x)$, $x\to p(x)$ and $x\to\dim H_j(x)$ are measurable. Nevertheless, if $\xi$ is an ergodic invariant measure for $f$, there is a subset $\Lambda(\xi)\subseteq\Lambda(f)$, such that $\xi(\Lambda(\xi))=1$ and there exist an integer $ p(\xi)\leq \dim M$, subspaces $H_1(\xi),...H_p(\xi)$, numbers $\chi_1(\xi)<\dots<\chi_p(\xi)$ such that for every $x\in\Lambda(\xi)$, we have
\begin{itemize}
\item $p(x)=p(\xi)$;
\item $\dim H_j(x)=\dim H_j(\xi)$, for every $j=1,\dots,p(\xi)$;
\item $\chi_j(x)=\chi_j(\xi)$, for every $j=1,\dots,p(\xi)$;
\end{itemize}

An ergodic measure $\xi$ is  {\em hyperbolic} if $\chi_j(\xi)\ne 0$, $j=1,\dots,p(\xi)$. In such case, for each $x\in\Lambda(\xi)$ we set
$$H^s(x)=\bigoplus_{\chi_j(\xi)<0}H_j(x), \mbox{ and }$$
$$H^u(x)=\bigoplus_{\chi_j(\xi)>0}H_j(x).$$
We have $\dim H^s(x)=s(\xi)$, $\dim H^u(x)=u(\xi)$ are constant and $s(\xi)+u(\xi)=\dim M$. The function $x\to H^s(x)$ and $x\to H^u(x)$ are measurables.
If $f$ is $C^r$, $r>1$, Pesin's Theory \cite{FHY81,Pe76,Pe77,PuSh89} guarantee the existence of invariant sub-manifolds $W^s(x)$, $W^u(x)$ tangent to $H^s(x)$  and $H^u(x)$ respectively.  More precisely, for every $x\in\Lambda(\xi)$ there is a $C^r$ embedded disk $W^s_{{\rm loc}}(x)$ through $x$ such that 
\begin{itemize}
\item $W^s_{{\rm loc}}(x)$ is tangent to $H^s(x)$ at $x$,
\item $f(W^s_{{\rm loc}}(x))\subseteq W^s_{{\rm loc}}(f(x))$, 
\item The stable set $W^s(x)=\cup_{n=0}^\infty f^{-n}(W^s_{{\rm loc}}(f^n(x))).$
\item There exist constant $C(x)>0$, $\tau(x)$ such that, for every $x_1,x_2\in W^s_{{\rm loc}}(x)$
\begin{equation}\label{eq:pesinsta}
{\rm dist}(f^k(x_1),f^k(x_2))\leq C(x)e^{-k\tau(x)}{\rm dist}(x_1,x_2).
\end{equation} 
\end{itemize}
 
The $C^r$ disk $W^s_{{\rm loc}}(x)$ is called {\em Pesin stable manifold}.  Similarly, every $x\in\Lambda(\xi)$ has an {\em Pesin unstable manifold} $W^u_{{\rm loc}}(x)$ satisfying the corresponding  properties with $f^{-1}$ in place of $f$.  

The Pesin manifolds above may be arbitrarily small, and they vary measurably on $x$. For any integer $n\geq 1$, we may find  {\em hyperbolic blocks} $\Lambda_n(\xi)\subseteq \Lambda(\xi)$ such that
\begin{itemize}
\item $\Lambda_n(\xi)\subseteq \Lambda_{n+1}(\xi)$,
\item $\xi(\Lambda_n(\xi))\to 1$, as $n\to \infty$.
\item   The the size of the embedded disk $W^s_{{\rm loc}}(x)$ is uniformly bounded from zero for each $x\in \Lambda_n(\xi)$. Moreover, for every $x\in \Lambda_n(\xi)$,  $C(x)<n$ and $\tau(x)>1/n$ in \eqref{eq:pesinsta}. Analogous properties are satisfied by the unstable Pesin's manifold $W^u_{{\rm loc}}(x)$.
\item The disk $W^s_{{\rm loc}}(x)$ and $W^u_{{\rm loc}}(x)$ vary continuously with $x\in\Lambda_n(\xi)$.
\end{itemize}

Most important, the holonomy maps associated to the Pesin stable lamination $\mathcal{W}^s_P=\{W^s_{{\rm loc}}(x)\}$ are absolutely continuous. More precisely, fix an integer $n\geq 1$, a hyperbolic block $\Lambda_n(\xi)$ and a point $x\in \Lambda_n(\xi)$. For $x_1$, $x_2\in W^s_{{\rm loc}}(x)$ close to $x$, let $\Sigma_1$ and $\Sigma_2$ be small smooth discs transverse to $W^s_{{\rm loc}}(x)$ at $x_1$ and $x_2$ respectively. The holonomy map 
$$\pi^s: \tilde{\Sigma_1}\subseteq \Sigma_1\to \Sigma_2$$
defined on the points $y_1\in\tilde{\Sigma_1}=\Sigma_1\cap \Lambda_n(\xi)$ by associate $\pi^s(y_1)$, the unique point in $\Sigma_2\cap W^s(y_1)$. If $f$ is $C^r$, $r>1$, then every holonomy map $\pi^s$ as before is absolutely continuous \cite{Pe76, PuSh89}. Of course, a dual statement holds for the unstable lamination.

In our setting, $f$ is a $C^r$-partially hyperbolic diffeomorphism, $r>1$, with  splitting $TM=E^{ss}\oplus E^c\oplus E^{uu}$, where $\dim E^\si=1$, $\si=ss,c,uu$. Let $\xi$ be an ergodic measure and we consider any point $x\in\Lambda(\xi)$. Then $p(x)=p(\xi)=3$  and $H_1(x)=E^s(x)$, $H_2(x)=E^c(x)$ and $H_3(x)=E^{uu}(x)$. Moreover
 $\chi_1(\xi)=:\lambda^s$, $\chi_3(\xi)=:\lambda^u$
 %\footnote{Comparar con la definicion y cotas que ponga Raul en sus previos} 
 and
 
$$\chi_2(x)=\lim_{n\to\pm\infty}\frac1n\log\|Df^{n}(x)|E^{c}_x\|=:\lambda^c(x),$$
is called the {\em center Lyapunov exponent} at $x$. If we take $x\in \mathcal{R}(\xi)\cap\Lambda(\xi)$,  since $\dim E^c=1$ we obtain that
\begin{equation}\label{eq:Leint}
\chi_2(\xi)=\lambda^c(\xi):=\lambda^c(x)=\int\log\|Df|E^c\|d\xi.
\end{equation}

If we assume $\lambda^c(\xi)<0$, then  $H^s(x)=E^{ss}(x)\oplus E^c(x)$ and $H^u(x)=E^{uu}(x)$ for every $x\in\Lambda(\xi)$. The local strong stable manifold  $W^{ss}_{{\rm loc}}(x)$ is an embedded curve inside the Pesin stable manifold $W^s_{{\rm loc}}(x)$ which is a surface. On the other hand, the Pesin unstable manifold $W^u_{{\rm loc}}(x)$ coincides with the strong unstable manifold $W^{uu}_{{\rm loc}}(x)$, for every $x\in\Lambda(\xi)$. Of course, analogous statement holds if we assume $\lambda^c(\xi)>0$.

Assume now that $f$ is partially hyperbolic and $\dim E^{uu}\geq 1$. An $f$-invariant probability measure $\mu$ is a {\em $u$-measure} if the conditional measures of $\mu$ with respect to the partition into local strong-unstable manifolds are absolutely continuous with respect to the Lebesgue measure along the corresponding local strong-unstable manifold. If $f$ is a $C^r$ partially hyperbolic diffemorphism, $r\geq2$, then there exist $u$-measures for $f$ \cite{PS}. Several properties of the $u$-measures are well know (see for instance \cite{BDV}, Section 11.2.3 and the references therein, for a detailed presentation of such properties). For instance, the support of any $u$-measure is a $u$-saturated, $f$-invariant, compact set. If $\mu$ is a $u$-measure, then its ergodic components are $u$-measures as well. Furthermore, the set of $u$-measures for $f$ is a compact, convex subset of the invariant measures. Moreover,  every physical measure for $f$ must be a $u$-measure.

It is well know that if $\mu$ is an ergodic $u$-measure with negative center Lyapunov exponent, then, $\mu$ is a physical measure \cite{Young:2002vc}. Conversely, if $\mu$ is a physical measure with negative center Lyapunov exponent, then $\mu$ is an ergodic $u$-measure.

\section{Examples}\label{sec:Examples}

In this section we show some examples that motivated this paper.  In the first example (Anosov times Morse-Smale)  there are no intermingled basins but there is a $u$-saturated set in the boundary of  one of them. Of course, we know a priori that this set consists of tori and it is not difficult to show that this situation is not robust. This example jointly with Kan's was a source of inspiration to obtain Theorem \ref{mTeo:A}. This is the easiest case where the theorem works. Observe that there is only one physical measure. 
In the second case (Kan-like example) the basins are intermingled.

\subsection{Toy Example}\label{ssec:Toy}
 In the 3-torus $\mathbb{T}^2\times \mathbb{S}^1$, we consider the $C^r$-diffeomorphism, $r\geq 2$, $F:\mathbb{T}^3\to\mathbb{T}^3$ defined by
$$F(x,t)=(Ax,\xi(t)),$$
where $A:\mathbb{T}^2\to \mathbb{T}^2$ is a linear Anosov diffeomorphism with eigenvalues $|\lambda_A^s|<1<|\lambda_A^u|$, and $\xi:\mathbb{S}^1\to \mathbb{S}^1$ is a Morse-Smale diffeomorphisms  with having exactly two hyperbolic fixed points, a source $p\in\mathbb{S}^1$ and a sink $q\in\mathbb{S}^1$,  satisfying $\partial W^u(p,\xi)=\{q\}$ and  $\partial W^s(q,\xi)=\{p\}$.  We assume that $F$ satisfies:
$$|\lambda_A^s|<|D\xi(t)|<|\lambda_A^u|, \mbox{ for every }t\in \mathbb{S}^1.$$
That means, $F$ is a partially hyperbolic diffeomorphism  exhibiting a center foliation by  compact leaves (circles).  Furthermore, $F$ has a foliation by smooth 2-tori tangent to the $E^s\oplus E^u$-sub-bundle. In particular, one of such leaves, the torus  $\mathbb{T}^2\times\{q\}$, is the only attractor of $F$. The  dynamics restricted to $\mathbb{T}^2\times\{q\}$ is hyperbolic, in fact, is given by $A$.  Then, it supports the unique hyperbolic $u$-measure $\mu_q$ for $F$ (actually the Lebesgue measure on  $\mathbb{T}^2\times\{q\}$) having negative center Lyapunov exponent and so,	 it is physical. If $\mathcal{B}_A(\mu_1)$ denotes the basin of $\mu_1$ in the 2-torus $\mathbb{T}^2\times\{q\}$ under the Anosov dynamics given by $A$ then, the basin of $\mu_q$  in $\mathbb{T}^3$  is
$$\mathcal{B}(\mu_q)=\mathcal{B}_A(\mu_q)\times\{\mathbb{S}^1\setminus\{p\}\},$$
which is an open set modulus a set of  zero Lebesgue measure in $\mathbb{T}^3$. The boundary of $\mathcal{B}(\mu_q)$ contains the invariant 2-torus $\mathbb{T}^2\times\{p\}$ which is the only hyperbolic repeller  of $F$. This invariant torus is also a $u$-saturated set, tangent to $E^s\oplus E^u$.
The  dynamics restricted to $\mathbb{T}^2\times\{p\}$ is again hyperbolic and then, it supports a $u$-measure $\mu_p$ for $F$ (actually Lebesgue measure on  $\mathbb{T}^2$) but it is not  physical. 

 Theorem A prevents the existence of such a $u$-saturated set from  being robust. After a typical $C^2$-perturbation, the new map $G$ is partially hyperbolic and dynamically coherent. In fact, $G$ has a  center foliation by compact leaves by classical results of normally hyperbolic foliations \cite{HPS}. 
% The plaque expansiveness  follows from to note that the quotient of $\mathbb{T}^3$ by the center leaves is the 2-torus $\mathbb{T}^2$ and the dynamics there is Anosov, which is expansive.  
% The dynamic of $G$ along the center foliation continues to be Morse-Smale, and the fixed points have their respective continuations. 

Typically $G$ does not preserve the invariant foliation by 2-tori tangent to $E^s\oplus E^u$, which exists for $F$. Nevertheless $G$ has two invariant compact  subset $\Lambda_p$ and $\Lambda_q$, the respective continuations of the hyperbolic basic sets $\mathbb{T}^2\times\{p\}$ and $\mathbb{T}^2\times\{q\}$. Of course, the dynamics of $F|\mathbb{T}^2\times\{p\}$ and $G|\Lambda_p$ are $C^0$-conjugated, so $\Lambda_p$ is (homeomorphic to) a continuous torus, and the dynamics of $G$ in $\Lambda_p$ is uniformly hyperbolic. The set $\Lambda_p$ remains to be a hyperbolic repeller and so  $s$-saturated, but  in general it is not  $u$-saturated.

Similar conclusions hold for $\Lambda_q$, the hyperbolic attractor of $G$. It is a topological 2-torus, $u$-saturated, and it supports the unique physical measure of $G$.  Note that the topological torus $\Lambda_p$ is contained in the boundary of the basin $\mathcal{B}(\mu_q^G)$, but, in general, $\Lambda_p$ is no longer a $u$-saturated set.

\subsection{Kan-like Examples}\label{ssec:Kan}

In \cite{Kan:1994kw} Kan  provided the first examples of  partially hyperbolic maps with intermingled basin. In this section we present the Kan's examples with some variations, following \cite{BDV}, Section 11.1.1.

\subsubsection{Kan's example: Endomorphism.}\label{sssec:Kan1}

\begin{figure}[h]
\includegraphics[scale=0.15]{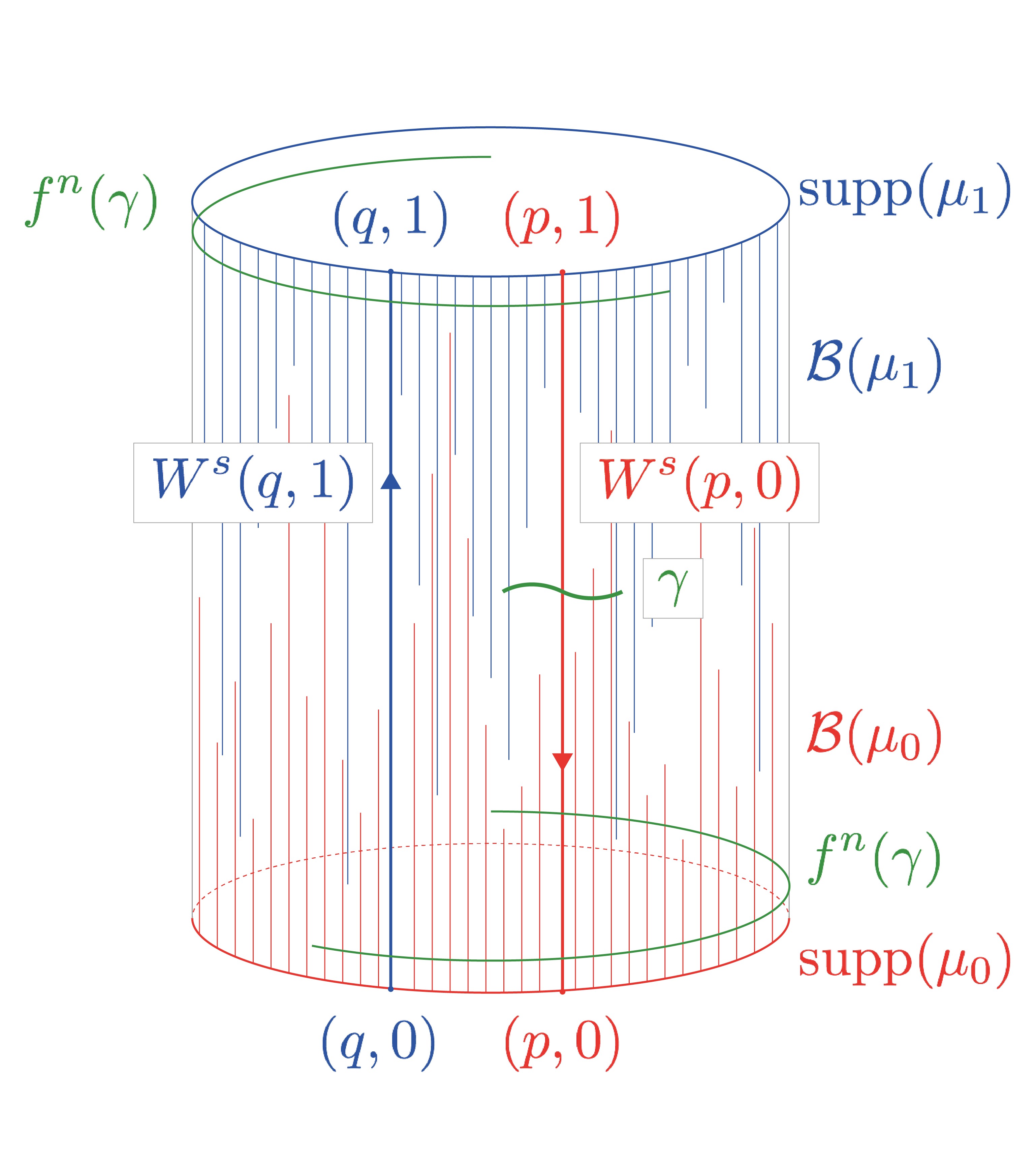} 
\caption{Kan example in the cylinder $\mathbb{S}^1\times [0,1]$}
\label{fig:kan}
\end{figure}

The Kan's example corresponds to a partially hyperbolic endomorphism defined on a surface with boundary exhibiting two intermingled hyperbolic physical measures.  Consider the cylinder $M=\mathbb{S}^1\times[0,1]$, and $K: M\to M$ the map defined by
$$K(\theta,t)=(k\theta \:({\rm mod}\: \mathbb{Z}),\varphi(\theta, t)),$$
where $k\geq 3$ is some integer,   $p, q\in \mathbb{S}^1$ are  two different fixed points of $\theta\to k\theta \:({\rm mod}\: \mathbb{Z})$  and $\varphi:M\to [0,1]$ is $C^r$, $r\geq 2$, satisfying the following conditions: 
\begin{enumerate}
\item[{[K1]}]For every $\theta\in \mathbb{S}^1$ we have $\varphi(\theta,0)=0$ and $\varphi(\theta,1)=1$. 
\item[{[K2]}]  The map $\varphi(p,\cdot):[0,1]\to[0,1]$ has exactly two fixed points, a hyperbolic source at $t=1$ and a hyperbolic sink in $t=0$. Analogously, the map  $\varphi(q,\cdot):[0,1]\to[0,1]$ has exactly two fixed points, a hyperbolic sink at $t=1$ and a hyperbolic source in $t=0$.
\item[{[K3]}]  For every $(\theta,t)\in M$, $|\partial_t\varphi(\theta,t))|<k$, and
\item[{[K4]}] $\displaystyle{\int\log |\partial_t \varphi (\theta,0)|\: d\theta<0}$ and  $\displaystyle{\int\log |\partial_t \varphi (\theta,1)|\: d\theta<0}$
\end{enumerate}
The dynamics along the $\theta$-direction is given by   $\theta \to k\theta \:({\rm mod}\: \mathbb{Z})$, so it is uniformly expanding. From [K3] we conclude that the map $K$ is partially hyperbolic: The derivative in the $t$-direction is dominated by the derivative in the $\theta$-direction. Condition [K1] means $K$ preserves the boundary.   Then, each one of the boundary circles  $\mathbb{S}^1\times\{0\}$  and $\mathbb{S}^1\times\{1\}$ supports an absolutely continuous invariant probability measure $\mu_0$ and $\mu_1$, respectively. Condition [K4] implies that $\mu_0$ and $\mu_1$ have negative Lyapunov exponent in the $t$-direction. So they are physical measures. Moreover, their basin are intermingled. Magic comes from condition [K2]: Take any curve $\gamma$ inside the open cylinder and transverse to the $t$-direction. We can assume, up to taking some forward iterates, that  $\gamma$ crosses (transversally) the segments $W^s(p,0)=\{p\}\times [0,1)$ and $W^s(q,0)=\{q\}\times (0,1]$. This is possible since $f$ is uniformly expanding along the $\theta$ direction and the angle between  $\gamma$ and the $t$-direction goes to $\frac\pi2$ due to the domination. Then, there is a forward iterate of $\gamma$ that intersects the basin of $\mu_0$, in a set of positive Lebesgue measure (in $\gamma$), because $\gamma$ intersects transversally $W^s(p,0)=W^{ss}(p)\times [0,1)$. Since $\gamma$ also intersects transversally $W^s(q,0)=W^{ss}(q)\times (0,1]$, then $\gamma$  intersects the basin of $\mu_1$ in a set of positive Lebesgue measure (See Figure~\ref{fig:kan}). Fubini's theorem completes the argument.

This example is robust among the maps defined on the cylinder preserving the boundaries. Indeed, for $r\geq 1$, any map $\tilde{K}: M\to M$, $C^r$ close to $K$ and preserving the boundaries can be written as
$$\tilde{K}(\theta,t)=(F(\theta, t),\psi(\theta, t)), $$
where $F:M\to \mathbb{S}^1$ is expanding along the $\theta$-direction and $\psi:M\to [0,1]$  preservs the boundaries, that means $\psi$ satisfies [K1]. Moreover, if $\psi$ is chosen $C^r$ close enough of $\varphi$, then also their derivatives $\partial_t\psi(\theta,t)$ and $\partial_t\varphi(\theta,t)$ are close for every $(\theta,t)\in M$ and so $\psi$ satisfies [K3] and [K4] above. The two different fixed points of $\theta\to k\theta \:({\rm mod}\: \mathbb{Z})$, $p, q\in \mathbb{S}^1$, have  continuations  $\tilde{p}, \tilde{q}\in \mathbb{S}^1$ and the map $\psi(\tilde{p},\cdot):[0,1]\to[0,1]$ has exactly two fixed points, a hyperbolic source at $t=1$ and a hyperbolic sink in $t=0$. Analogously, the map  $\psi(\tilde{q},\cdot):[0,1]\to[0,1]$ has exactly two fixed points, a hyperbolic sink at $t=1$ and a hyperbolic source in $t=0$. Then, $\psi$ satisfies [K2]. Arguing as before, we conclude that $\tilde{K}$ exhibits two intermingled hyperbolic physical measures supported on the boundary.

\subsubsection{Kan's example:  Diffeomorphisms on a manifold with boundary}\label{sssec:Kan2}

The next example, corresponds to a partially hyperbolic diffeomorphism defined on a 3-manifold with boundary exhibiting two intermingled physical measures.  The idea is to adapt the previous example,   replacing $\mathbb{S}^1$ with the torus $\mathbb{T}^2$ and the expanding map $\theta\to k\theta \:({\rm mod}\: \mathbb{Z})$ with a hyperbolic automorphism of the 2-torus having at least two fixed points.  More precisely, we can consider $N=\mathbb{T}^2\times[0,1]$ and diffeomorphisms
$$K_D(z,t)=(Az,\psi(\theta, t)),$$
where $A:\mathbb{T}^2\to\mathbb{T}^2$ is a hyperbolic automorphism, and $\psi:N\to [0,1]$ is $C^r$, $r\geq 2$, satisfying the following conditions: 
\begin{enumerate}
\item[{[KD1]}] For every $z\in \mathbb{T}^2$ we have $\psi(z,0)=0$ and $\psi(z,1)=1$. 
\item[{[KD2]}]  For $p,q\in\mathbb{T}^2$, fixed points of $A$, we assume that the map $\psi(p,\cdot):[0,1]\to[0,1]$ has exactly two fixed points, a source at $t=1$ and a sink in $t=0$. Analogously, the map  $\psi(q,\cdot):[0,1]\to[0,1]$ has exactly two fixed points, a sink at $t=1$ and a source in $t=0$.
\item[{[KD3]}]  For every $(z,t)\in M$, $\|A^{-1}\|^{-1}<|\partial_t\psi(x,t))|<\|A\|$, and
\item[{[KD4]}] $\displaystyle{\int_{\mathbb{T}^2}\log |\partial_t \psi (z,0)|\: dz<0}$ and  $\displaystyle{\int_{\mathbb{T}^2}\log |\partial_t \psi (z,1)|\: dx<0}$
\end{enumerate}
As before, the dynamics along the $z$-direction of $K_D$ is   uniformly hyperbolic.  From [KD3] we conclude that the map $K_D$ is partially hyperbolic: The derivative in the $t$-direction is dominated by the derivative in the unstable direction of $A$ and the stable direction of $A$ is dominated by the derivative in the $t$-direction. Condition [KD1] means $K_D$ preserves each boundary torus.   Then both boundary torus  $\mathbb{T}^2\times\{0\}$  and $\mathbb{T}^2\times\{1\}$ support the measures  $\mu_0$ and $\mu_1$ corresponding to the Lebesgue measure in the torus. Condition [KD4] implies that $\mu_0$ and $\mu_1$ have negative Lyapunov exponent in the center direction. So they are physical measures.

As before, their basins are intermingled. The argument is the same: Take any curve $\gamma$ in the interior of $N$ and transverse to the 
$E^{cs}$ distribution. Up to some forward iterates,  $\gamma$ crosses (transversally) the surfaces $W^s_{loc}(p,0)=W^{ss}_{loc}(p)\times [0,1)$ and $W^s_{loc}(q,0)=W^{ss}_{loc}(q)\times(0,1]$. This is possible since $f$ is uniformly expanding along the unstable direction and the domination improves the angle between  $\gamma$  and  the center-stable direction. Then, there is a forward iterate of $\gamma$ that intersects the basin of $\mu_0$ in a set of positive Lebesgue measure (in $\gamma$), because $\gamma$ intersects transversally the  stable manifold $W_{loc}^s(p,0)$. Since $\gamma$ also intersects transversally the  stable manifold $W^s(q,0)$, then $\gamma$  intersects the basin of $\mu_1$ in a set of positive Lebesgue measure. Fubini's theorem complete the argument.

As before, this example is robust among the diffeomorphisms defined on $N$ preserving the boundary tori. 

\subsubsection{Kan-like example: Diffeomorphisms on a boundaryless manifold}\label{sssec:Kan3}

The same construction can be done if $N$ is replaced with $\mathbb{T}^3=\mathbb{T}^2\times \mathbb{S}^1$ (or even the mapping torus of a hyperbolic diffeomorphism) and $\psi:N \to [0,1]$ is replaced with $\varphi:\mathbb{T}^2\times \mathbb{S}^1\to \mathbb{S}^1$. Then, the four conditions are: 

\begin{enumerate}
\item[{[KB1]}] For every $z\in \mathbb{T}^2$ we have $\varphi(z,0)=0$ and $\varphi(z,\frac12)=\frac12$. 
\item[{[KB2]}]  For $p,q\in\mathbb{T}^2$, fixed point of $A$, we assume that the map $\varphi(p,\cdot):\mathbb{S}^1\to\mathbb{S}^1$ has exactly two fixed points, a source at $t=\frac12$ and a sink in $t=0$. Analogously, the map  $\varphi(q,\cdot):\mathbb{S}^1\to\mathbb{S}^1$ has exactly two fixed points, a sink at $t=\frac12$ and a source in $t=0$.
\item[{[KB3]}]  For every $(z,t)\in M$, $\|A^{-1}\|^{-1}<|\partial_t\varphi(x,t))|<\|A\|$, and
\item[{[KB4]}] $\displaystyle{\int_{\mathbb{T}^2}\log |\partial_t \varphi (z,0)|\: dz<0}$ and  $\displaystyle{\int_{\mathbb{T}^2}\log |\partial_t \varphi (z,\frac12)|\: dx<0}$.
\end{enumerate}

Exactly the same proof gives that the basins of the Lebesgue measures of the boundary tori are intermingled. The difference is that this phenomenon is no longer robust. In fact, there exists a unique physical measure after most perturbations (see, for instance, \cite{DVY}).

Recently, Bonatti and Potrie announced that they are able to construct diffeomorphisms on the torus $\mathbb{T}^3$ with exactly $k\geq 2$ hyperbolic physical measures $\mu_1,\dots,\mu_k$ whose basins are all intermingled (and dense on the whole torus), in fact, for every open set $A\subseteq \mathbb{T}^3$ and every $i\ne j\in\{1,\dots, k\}$

$${\rm Leb}(A\cap \mathcal{B}(\mu_i))>0\qquad \mbox{ and }\qquad {\rm Leb}(A\cap \mathcal{B}(\mu_j))>0.$$

Their example is partially hyperbolic in the following broad sense: the tangent space has an invariant splitting $T\mathbb{T}^3=E^{cs}\oplus E^u$ where $E^u$ dominates $E^{cs}$ but the sub-bundle $E^{cs}$ is indecomposable into dominated sub-bundles.

We  remark that  partially hyperbolic diffeomorphisms  on surfaces do not admit intermingled hyperbolic physical measures \cite{Hertz:2011vu}. The situation is different in the absence of domination as showed by Fayad \cite{F2003}. Inspired in the Fayad example, Melbourne and Windsor \cite{MELBOURNE:2005hn} give a family of $C^\infty$-diffeomorphisms on $\mathbb{T}^2\times \mathbb{S}^2$ with arbitrary number of physical measures with intermingled basins.

Motivated by the latter situation, we say that a partially hyperbolic diffeomorphism $f$ is a {\em Kan-like differmorphisms} if there exist, at least, two hyperbolic physical measures with intermingled basins.

\section{Proof of Theorem~\ref{mTeo:A} and Corollary~\ref{mcor:B}}\label{sec:proofs}

%For every $x\in M$, the set
%$$P(x)=W^s(W^u(x)):=\bigcup_{y\in W^u(x)}W^s(y),$$
%is a topological surface in the universal covering, and it is topologically transversal to the center direction $E^c$\cite{H2009,H2013}.

Let $f\in{\rm Diff}^r(M)$, $r\geq 2$, be partially hyperbolic and dynamically coherent with compact center leaves. Let $\mu$ be a hyperbolic physical measure for $f$ with  $\lambda^c(\mu)<0$. For further use let $\Lambda=\cup_n\Lambda_n$ where $\Lambda_n$ are Pesin blocks and $\mu(\Lambda)=1$. We assume that $\Lambda$  is invariant and its points are regular both in the sense of Pesin's Theory as in the sense of Birkhoff's Theorem. Moreover, we will assume that every $x\in \Lambda_n$ is a Lebesgue density point of $W^u(x)\cap \Lambda_n$.

For $E\subseteq M$ measurable, $W^s(E)$ denotes the union of Pesin's stable manifolds $W^s(x)$ of  points  $x\in E$. Observe that $W^s(E)$ is invariant if $E$ is invariant.

First, for the sake of completeness, we will prove the following lemma. We thank the referee for provide us the argument of the proof.

\begin{lema}\label{l1} 
$\mathcal{B(\mu)}\subset \overline{W^s(\Lambda)}$. 
\end{lema}

\proof Suppose $x\in \mathcal{B(\mu)}$. Fix $m \geq 1$. Then, it is not difficult to see there is a sequence $n_k$ such that the distance between $f^{n_k}(x)$ and $\Lambda_m$ converges to $0$. Indeed, if there is $\delta>0$ such that the distance  of  $f^n(x)$ to $\Lambda_m$ is greater than $0$ you can construct  a continuous which that takes the value $1$ for every point in $\Lambda_m$ and $0$ if the distance to $\Lambda_m$ is greater or equal to $\delta$. Since $\Lambda_m$ has positive $\mu$-measure this contradicts the fact that $x\in \mathcal{B(\mu)}$.

As Pesin stable manifolds are of uniform size for points in $\Lambda_m$, there is
$y_k\in W^u_{loc}(f^{n_k})\cap W^s(\Lambda)$ for any $k$ large enough. Clearly $f^{-n_k}(y_k)$ converges to $x$ proving the lemma.

\endproof

Denote by $M_c$  the space of center curves, that is, the quotient space obtained by the relation of equivalence $y\sim x$ if they are in the same center manifold. We denote by  $X$  the space of compact subsets of $M$. Given a $u$-saturated closed subset $K\subseteq M$,  we define  the function $\Phi_K: M_c\to X$  by $\Phi_K(\bar x)=K\cap \bar x$.  Observe that this intersection is nonempty for every $\bar x\in M_c$.

Since $K$ is closed we have that $\Phi_K$ is an upper semicontinuous function. This implies that $\Phi_K$ has a residual set of points of continuity. 

On the other hand, since $K$ is saturated by strong unstable leaves and the unstable holonomy is continuous,   the set of continuity points of $\Phi_K$ is also saturated by strong unstable leaves. More precisely, if $\bar x$ is a  point of continuity of $\Phi_K$, then for every $y\in W^u(\bar{x})$ we have that   $\bar{y} \in M_c$ is also a point of continuity of $\Phi_K$.

%We define a function $\Psi: M_c \to X$ given by $\Psi(\bar x) = \supp(\mu)\cap \bar x$. Since $\supp(\mu)$ is closed we have that $\Psi$ is an upper semicontinuous function. This implies that it has a residual set of points of continuity. 
%
%Observe that the facts that $\supp(\mu)$ is saturated by strong unstable leaves and the unstable holonomy is continuous imply the the set of continuity points of $\Psi$ is also saturated by strong unstable leaves. 

\begin{lema}\label{l2}
For every $x\in W^s(\Lambda)$, there is a center arc $[x,y]_c\subseteq W^s(y)$ with $y\in{\rm supp}(\mu)$.
\end{lema}

\proof

Let $x\in W^s(\Lambda_m)$. Taking iterates for the future, and recalling that almost every point returns infinitely many times to a positive measure set,  we can assume that $x\in W^s_\e(y')$ with $y'\in \Lambda_m$  where $\e$ is the uniform size of the Pesin stable manifolds of the points of the block $\La_m$.
Close to $y'$ we take $z\in W^s(y')\cap \supp(\mu)$, with ${\rm dist}(y',z)<\frac\e{10}$, and such that $\bar z=W^c(z)$ is a continuity point of $\Phi_{\supp(\mu)}$. In particular, there is a $\de>0$  such that, if ${\rm dist}(z,w)<\delta$ then, there exists $p\in W^c(w)\cap \supp(\mu)$ with ${\rm dist}(p,z)<\frac\e{10}$.

Let $H=\La_m \cap W^{uu}_{\frac\e{10}}(y')$ and $G=W^s_\epsilon(H)\cap B_{\de/2}(z)\cap \supp(\mu)$. The absolute continuity of the partition by Pesin' stable manifolds implies that $\mu(G)>0$.
Then, the ergodicity of the measure implies that there are infinitely many iterates of $y'$ that belong to $G$. In particular, there is an $n$ such that $f^n(y')\in G$ and ${\rm dist}(f^n(x), f^n(y'))<\de/2$. Thus, we obtain that $f^n(x)\in B_\de(z)\cap W^s_\e(H)$.
The fact that $f^n(x)\in B_\de(z)$ implies that there is $v\in W^c(f^n(x))\cap\supp\mu$, such that ${\rm dist}(f^n(x),v)< \frac\e{10}$. Since $f^n(x) \in W^s_\e(H)$ we have that corresponding center arc $[f^n(x),v]_c$ is completely contained in a Pesin stable manifold. We take $y=f^{-n}(v)$ and this gives the conclusion of the lemma for the points of $W^s(\La)$. 
\endproof

In what follows we consider $K\subseteq \mathbb{T}^3$ satisfying the hypotheses in Theorem~\ref{mTeo:A}. That is, $K$ is  a compact,  $f$-invariant and $u$-saturated subset such that $K\subseteq \overline{\mathcal{B(\mu)}}\setminus\supp\mu$. Our strategy to prove Theorem \ref{mTeo:A} will be to study the intersections of the set $K$ with the center manifolds of $f$.  

\begin{lema}\label{lemah} There is a $h>0$ such that if we have three distinct points $x, y, z \in \Phi_K(\bar w)$ then at least two of them are a $c$-distance larger than $h$.
\end{lema}

\proof
As we have already mentioned we will use that $W^{ss}(W^{uu}(x))$, when considered in the universal cover, is topological surface topologically transverse to the center leaves \cite{H2009,H2013}.

\begin{figure}[h]
\includegraphics[scale=0.2]{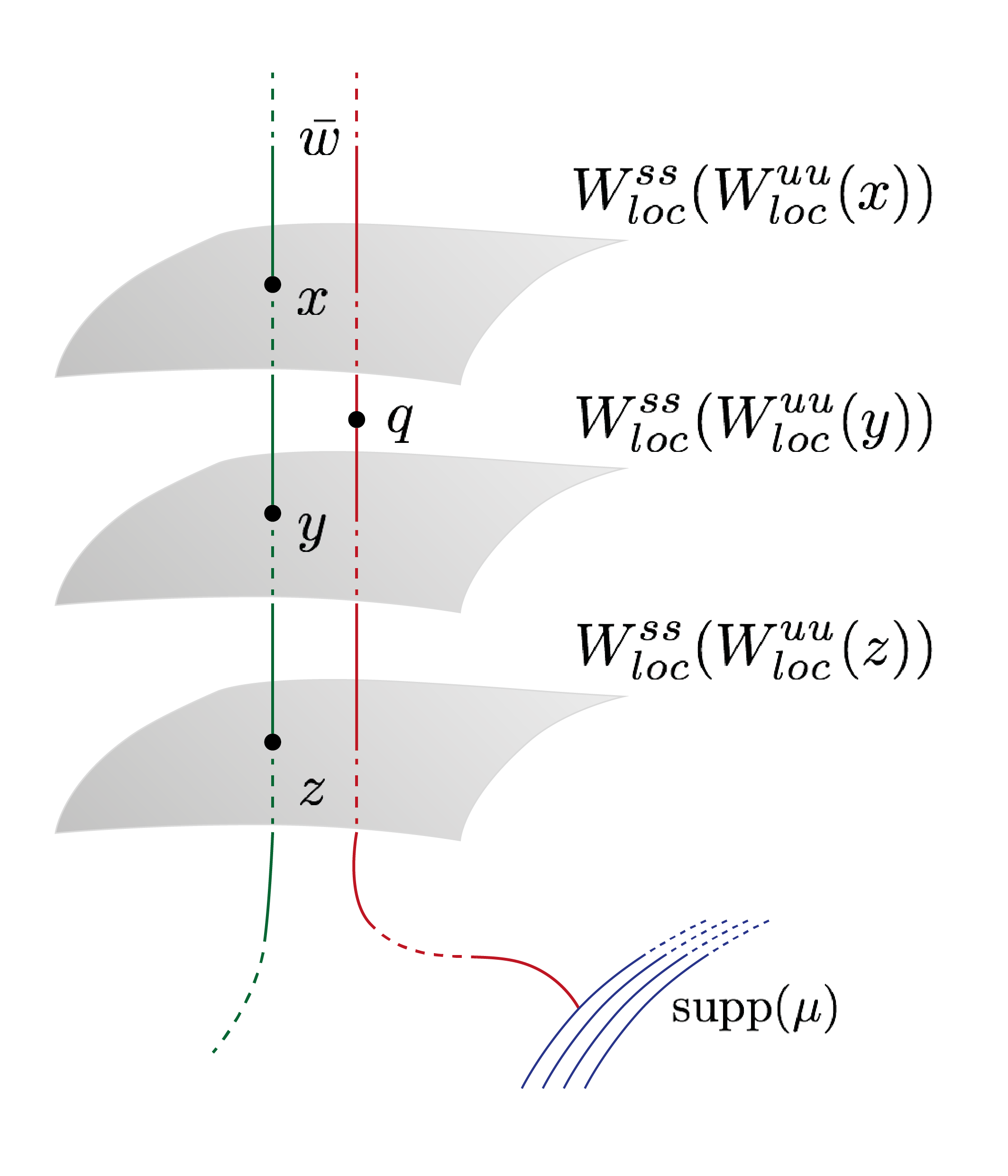} 
\caption{The center arc (in red) joining $q$ and $\supp(\mu)$ must intersect the local planes.}
\label{fig:1}
\end{figure}

Let's begin with the proof.  Suppose on the contrary that for every $h$ there are $w$ and  three points $x, y, z \in \Phi_K(\bar w)$ with ${\rm dist}_c(u,v)<h$ for every pair of points  $\{u,v\}\subset\{x,y, z\}$. Take the topological surfaces $W_{loc}^{ss}(W_{loc}^{uu}(x))$,  $W_{loc}^{ss}(W_{loc}^{uu}(y))$ and $W_{loc}^{ss}(W_{loc}^{uu}(z))$. Without loose of generality we can assume that $y$ is in the center arc that joins $x$ and $z$ and has length less than $h$. Take $k>0$ such that ${\rm dist}(K, \supp\mu)>k$ and suppose that $h\ll k$.  Since $y\in \overline{\B(\mu)} $, Lemma \ref{l1} implies that  it can be approximated by a point $q$ belonging to $W^s(\Lambda)$. By Lemma \ref{l2} we have that $q$ can be joined to $\supp(\mu)$ by a center arc completely contained in $\B(\mu)$. Observe that $q$ is very close to $y\in K$ and then, the length of this center arc is greater than, say, $k/2$. Still much larger than $h$. This implies that the center arc joining $q$ and $\supp(\mu)$ must intersect either $W_{loc}^{ss}(W_{loc}^{uu}(x))$ or $W_{loc}^{ss}(W_{loc}^{uu}(z))$ (See Figure~\ref{fig:1}).

This is a contradiction because these sets are in the complement of $\B(\mu)$, since the $\omega$-limits of its points are in $K$. This ends the proof of the lemma.

\endproof

The preceding lemma has an immediate and important consequence that we state as proposition.

\begin{prop}\label{finito} Let $\mu$ be an ergodic  $u$-measure with negative center exponent and $K$ an invariant $u$-saturated set such that $K\subseteq \overline{\B(\mu)}\setminus \supp\mu$. Then, the intersection of $K$ with each center manifold consists of finitely many points.
\end{prop}

Our next  lemma says that the number of points of the intersection of a $u$-minimal subset $J$ of $K$ with each center manifold is constant. 
 
\begin{lema}\label{constant} Let $J\subseteq K$ be $u$-minimal compact set. Then $\#\Phi_J(\bar x)$ does not depend on $\bar x$
\end{lema}

\proof
We want to show that the function $\#\Phi_J$ is constant in an open set. If this is the case, the $u$-minimality of $J$ and the $u$-invariance of $\#\Phi_J$ will imply the proposition. 

Observe that, a priori, the semicontinuity of $\Phi_J$ does not imply directly the proposition because it is not enough to conclude the semicontinuity of $\#\Phi_J$.

 Let $\bar x$ be a point of continuity of $\Phi_J$. Continuity at $\bar x$ implies that $\#\Phi_J(\bar y)\geq \#\Phi_J(\bar x)$ if $\bar y$ is close enough to $\bar x$. The $u$-minimality, again,  implies the inequality for every $\bar y\in M_c$. 
 
 Suppose that the function $\#\Phi_J$ is not constant. Then, there is a dense set $D\subseteq M_c$ such that for $\bar y\in D$ we have that $\#\Phi_J(\bar y)> \#\Phi_J(\bar x)$. Continuity at $\bar x$ implies that there are a point $x\in \bar x$, a sequence $\bar y_n \to \bar x$ and  for each integer $n\geq 1$, a pair of points $y_n^1$, $y_n^2\in \bar y_n \cap J$ so that both sequences $(y^i_n)$, $i=1,2$, converge to $x$. Then, taking $N$ large enough we can choose a center curve with two points $y_N^1:=y^1$ and $y_N^2:=y^2$ a very small $c$-distance. We will argue in a similar way to the arguments of the proof of Lemma \ref{lemah}. We want to obtain three points that are very close to each other in the same center manifold and surfaces through them that are not in $\B(\mu)$, to arrive to a contradiction with Lemma \ref{l2}.
 
\begin{figure}[h]
\includegraphics[scale=0.2]{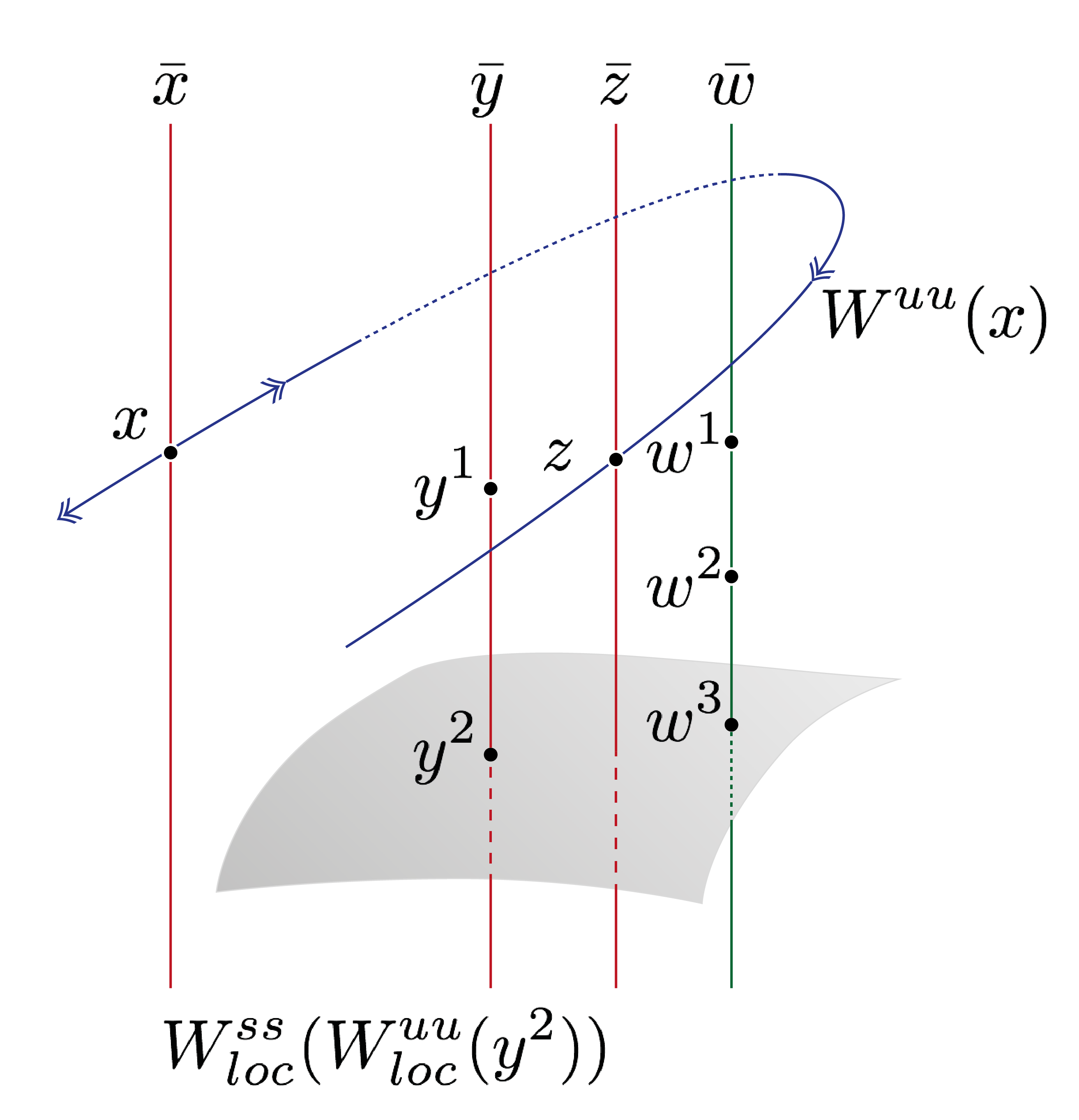} 
\caption{$W_{loc}^{ss}(W^{uu}_{loc}(y^2))$ intersects $\bar w$ in a point $w^3$ that is close to $w^1$ and $w^2$.}
\label{fig:3}
\end{figure}

 Since $J$ is $u$-minimal we can find $z\in W^{uu}(x)$ very close to $y^1$. Continuity of the holonomy gives that there are center manifolds converging to the center manifold of $z$ and pairs of points $w^1_n, \,w^2_n$ of $J$ in each of these center manifolds converging to $z$. Finally, fix an integer $L\geq 1$ large enough (in such a way that the $c$ distance between $w^1_L$ and $w^2_L$ is much smaller than the one between $y^1$ and $y^2$) and call $w^1=w^1_L$ and $w^2=w^2_L$. Denote  $\bar w$ the center leaf that contains $\{w^1,w^2\}$. Because of the choices we have made, $W_{loc}^{ss}(W^{uu}_{loc}(y^2))$ intersects $\bar w$ in a point $w^3$ that is close to $w^1$ and $w^2$ but at a greater distance than ${\rm dist}_c(w^1,w^2)$. That means that one of the two points $w^1,\, w^2$ lies in between the other two  (See figure \ref{fig:3}).
 
Now, arguing as in Lemma \ref{lemah} we arrive to a contradiction.

\endproof

\proof[Proof of Theorem \ref{mTeo:A}] Let $J\subseteq K$ be $u$-minimal and closed.  
Lemma \ref{constant} shows that $J$ is locally the graph of a continuous function and then, it is a closed topological surface topologically transverse to the center foliation. Since it is foliated by unstable leaves, that are lines, we have that $J$ is a torus. Moreover, Proposition \ref{finito} implies that the torus $J$ is periodic. Thus, all that remains is to prove that the strong stable manifolds of the points of $J$ are completely contained in $J$.

\begin{figure}[h]
\includegraphics[scale=0.2]{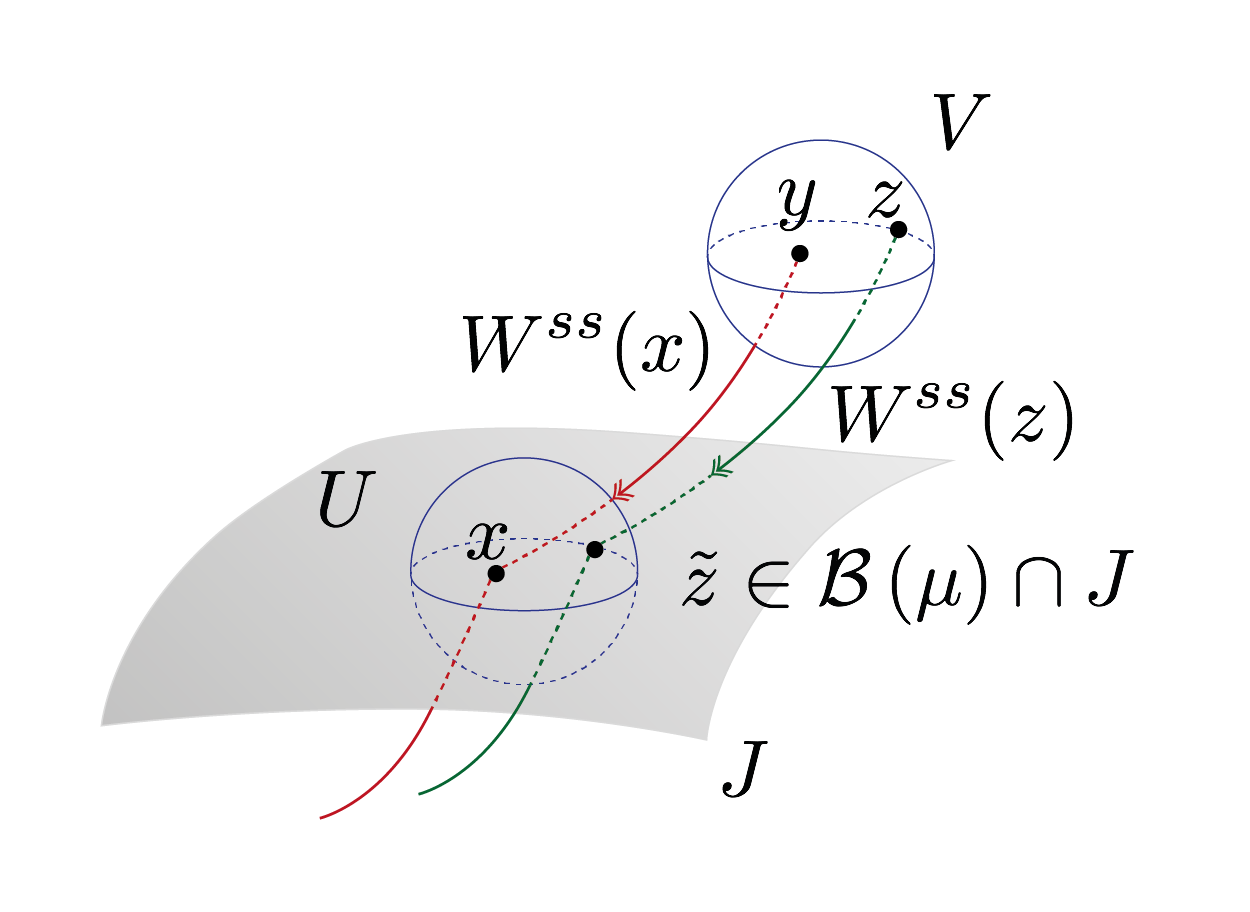} 
\caption{Graphic representation of the proof of Theorem \ref{mTeo:A}.}
\label{fig:2}
\end{figure}

As $J\subseteq K$ is periodic, we can take  an iterate $n\geq 1$ such that  $f^n(J)=J$. By simplicity we assume that $n=1$.  Suppose that there is a point $x\in J$  such that its strong stable manifold $W^{ss}(x)$ has a point $y$ that does not belong to $J$. Since $J$ is closed, there exists an open neighbourhood $V\subseteq M$ of $y$ such that $V\cap J=\emptyset$.  By the continuity of the strong stable foliation, reducing $V$ if necessary, we can find an open neighbourhood $U\subseteq M$  of $x$ with the property that the strong stable manifold of every point in $V$ has a point  in $U$, in particular, in $J$. We know that $J\subseteq K\subseteq \overline{\B(\mu)}\setminus\supp\mu$, then  $V\cap \B(\mu)\ne \emptyset$. Hence, there is $z\in V\cap \B(\mu)$ and if we take $\tilde{z}\in W^{ss}(z)\cap J$, then $\tilde{z}\in \B(\mu)$ (See Figure \ref{fig:2}). 

In particular, $\tilde{z}\in J$ and its omega limit is contained in $\supp(\mu)$.  %This means that  the sequence of measures $$\mu_n=\frac{1}{n}\sum_{j=0}^{n-1}\delta_{f^j(\tilde{z})}$$ converge to $\mu$ in the weak* topology. 
Since $J$ is $f$-invariant, then $\emptyset\ne\supp(\mu)\cap J\subseteq K$ 
which contradicts the hypotheses $J\subseteq K\subseteq\overline{\B(\mu)}\setminus\supp\mu$. This finishes the proof of the Theorem~\ref{mTeo:A}.

\endproof

\proof[Proof of Theorem \ref{mcor:B}]   Let  $\mu$ and $\nu$ be two hyperbolic physical measures.  Recall that  their supports are compact,  $f$-invariant and $u$-saturated subsets. 

First of all, observe that neither $\mu$ nor $\nu$ can have positive center Lyapunov exponent. This is a consequence of the well-known fact that under our hypotheses the basin of attraction of such a measure would be essentially open (See  for instance \cite{BDP} where the conservative case is discussed with details and recently \cite{AV2015} for a discussion about  the non conservative case.) . 

Suppose that the center exponents are negative. If their basins are intermingled then $\supp\:\nu\subseteq \overline{\mathcal{B(\mu)}}\setminus\supp\:\mu$. Indeed, it is not difficult to see that the definition of intermingled basins implies that there is a point of the stable manifold (in the sense of Pesin) of a regular point of $\nu$ that is accumulated by points of the basin of $\mu$. Since $\nu$ is ergodic the orbit of a regular point is dense in its support. By forward iteration we obtain the desired inclusion.  Then, as  consequence of Theorem~\ref{mTeo:A} applied to $K=\supp \nu$, $f$ is not accessible. As mentioned above accessibility is an open an dense property, and then we obtain the first assertion.  

For the second statement, the works of A. Hammerlindl \cite{H2009}  and R. Potrie \cite{PotJMD} proved that the center foliation of every dynamically coherent partially hyperbolic diffeomorphism on the 3-torus is homeomorphic to the corresponding foliation of a linear toral automorphism. As a consequence, there are two possibilities: either the center foliation is by circles or the diffeomorphism is homotopic to a hyperbolic automorphism, it is  always dynamically coherent  and the center foliation is by lines. We have already studied the first case. In the second case, 
Potrie \cite{Pot2015} (see also \cite{U}) proved that if $f$ is isotopic to a hyperbolic automorphism, 
there is a unique minimal $u$-saturated set. This implies that $f$ has at most one physical measure with negative center exponent. 

\endproof

\subsubsection*{Acknowledgement} The authors would like to thank the anonymous reviewer for their helpful and constructive comments.

\bibliographystyle{plain}

%\bibliography{URVA}

\begin{thebibliography}{10}

\bibitem{AV2015}
M.~{Andersson} and C.~H. {V{\'a}squez}.
\newblock {On mostly expanding diffeomorphisms}.
\newblock {\em ArXiv e-prints}:1512.01046, December 2015.

\bibitem{BP}
L.~Barreira and Y.~Pesin.
\newblock {\em Lyapunov exponents and smooth ergodic theory}, volume~23 of {\em
  University Lecture Series}.
\newblock American Mathematical Society, Providence, RI, 2002.

\bibitem{BDV}
C.~Bonatti, L.~D{\'i}az, and M.~Viana.
\newblock {\em Dynamics beyond uniform hyperbolicity}, volume 102 of {\em
  Encyclopaedia of Mathematical Sciences}.
\newblock Springer-Verlag, Berlin, 2005.
\newblock A global geometric and probabilistic perspective, Mathematical
  Physics, III.

\bibitem{BM2008}
A.~Bonifant and J.~Milnor.
\newblock Schwarzian derivatives and cylinder maps.
\newblock In {\em Holomorphic dynamics and renormalization}, volume~53 of {\em
  Fields Inst. Commun.}, pages 1--21. Amer. Math. Soc., Providence, RI, 2008.

\bibitem{BBI2009}
M.~Brin, D.~Burago, and S.~Ivanov.
\newblock Dynamical coherence of partially hyperbolic diffeomorphisms of the
  3-torus.
\newblock {\em J. Mod. Dyn.}, 3(1):1--11, 2009.

\bibitem{BP1974}
M.~Brin and Y.~Pesin.
\newblock Partially hyperbolic dynamical systems.
\newblock {\em Izv. Akad. Nauk SSSR Ser. Mat.}, 38:170--212, 1974.

\bibitem{BDP}
K.~Burns, D.~Dolgopyat, and Y.~Pesin.
\newblock Partial hyperbolicity, {L}yapunov exponents and stable ergodicity.
\newblock {\em J. Statist. Phys.}, 108(5-6):927--942, 2002.
\newblock Dedicated to David Ruelle and Yasha Sinai on the occasion of their
  65th birthdays.

\bibitem{BHHTU}
K.~Burns, F.~Rodriguez~Hertz, J.~Rodriguez~Hertz, A.~Talitskaya, and R.~Ures.
\newblock Density of accessibility for partially hyperbolic diffeomorphisms
  with one-dimensional center.
\newblock {\em Discrete Contin. Dyn. Syst.}, 22(1-2):75--88, 2008.

\bibitem{DVY}
D.~Dolgopyat, M.~Viana and J.~Yang.
\newblock Geometric and measure-theoretical structures of maps with mostly
  contracting center.
\newblock {\em Comm. Math. Phys.} 341(3): 991-1014, 2016.

\bibitem{FHY81}
A.~Fathi, M.~Herman, and J.~Yoccoz.
\newblock A proof of {P}esin's stable manifold theorem.
\newblock In {\em Geometric dynamics ({R}io de {J}aneiro, 1981)}, volume 1007
  of {\em Lecture Notes in Math.}, pages 177--215. Springer, Berlin, 1983.

\bibitem{F2003}
B.~Fayad.
\newblock Topologically mixing flows with pure point spectrum.
\newblock In {\em Dynamical systems. {P}art {II}}, Pubbl. Cent. Ric. Mat. Ennio
  Giorgi, pages 113--136. Scuola Norm. Sup., Pisa, 2003.

\bibitem{H2009}
A.~Hammerlindl.
\newblock {\em Leaf conjugacies on the torus}.
\newblock ProQuest LLC, Ann Arbor, MI, 2009.
\newblock Thesis (Ph.D.)--University of Toronto (Canada).

\bibitem{H2013}
A.~Hammerlindl.
\newblock Leaf conjugacies on the torus.
\newblock {\em Ergodic Theory Dynam. Systems}, 33(3):896--933, 2013.

\bibitem{HP2014}
A.~Hammerlindl and R.~Potrie.
\newblock Pointwise partial hyperbolicity in three-dimensional nilmanifolds.
\newblock {\em J. Lond. Math. Soc. (2)}, 89(3):853--875, 2014.

\bibitem{HPS}
M.~Hirsch, C.~Pugh, and M.~Shub.
\newblock {\em Invariant manifolds}.
\newblock Springer-Verlag, Berlin, 1977.
\newblock Lecture Notes in Mathematics, Vol. 583.

\bibitem{IKS08}
Y.~Ilyashenko, V.~Kleptsyn, and P.~Saltykov.
\newblock Openness of the set of boundary preserving maps of an annulus with
  intermingled attracting basins.
\newblock {\em J. Fixed Point Theory Appl.}, 3(2):449--463, 2008.

\bibitem{Kan:1994kw}
I.~Kan.
\newblock Open sets of diffeomorphisms having two attractors, each with an
  everywhere dense basin.
\newblock {\em Bull. Amer. Math. Soc. (N.S.)}, 31(1):68--74, 1994.

\bibitem{KS2011}
V.~Kleptsyn and P.~Saltykov.
\newblock On {$C^2$}-stable effects of intermingled basins of attractors in
  classes of boundary-preserving maps.
\newblock {\em Trans. Moscow Math. Soc.}, vol 17, pages 193--217, 2011.

\bibitem{M87}
R.~Ma{\~n}{\'e}.
\newblock {\em {Ergodic theory and differentiable dynamics}}, volume~8 of {\em
  Ergebnisse der Mathematik und ihrer Grenzgebiete (3) [Results in Mathematics
  and Related Areas (3)]}.
\newblock Springer-Verlag, Berlin, Berlin, Heidelberg, 1987.

\bibitem{MELBOURNE:2005hn}
I.~Melbourne and A.~Windsor.
\newblock A {$C^\infty$} diffeomorphism with infinitely many intermingled
  basins.
\newblock {\em Ergodic Theory Dynam. Systems}, 25(6):1951--1959, 2005.

\bibitem{NT2001}
V.~Ni{\c{t}}ic{\u{a}} and A.~T{\"o}r{\"o}k.
\newblock An open dense set of stably ergodic diffeomorphisms in a neighborhood
  of a non-ergodic one.
\newblock {\em Topology}, 40(2):259--278, 2001.

\bibitem{Okunev}
A.~{Okunev}.
\newblock {Milnor attractors of circle skew products}.
\newblock {\em ArXiv e-prints}:1508.02132., August 2015.

\bibitem{Oseledec:1968tk}
V.~Oseledec.
\newblock {A multiplicative ergodic theorem. Characteristic Ljapunov, exponents
  of dynamical systems}.
\newblock {\em Trudy Moskovskogo Matemati\v ceskogo Ob\v s\v cestva},
  19:179--210, 1968.

\bibitem{Pe76}
Y.~Pesin.
\newblock Families of invariant manifolds that correspond to nonzero
  characteristic exponents.
\newblock {\em Izv. Akad. Nauk SSSR Ser. Mat.}, 40(6):1332--1379, 1440, 1976.

\bibitem{Pe77}
Y.~Pesin.
\newblock Characteristic {L}japunov exponents, and smooth ergodic theory.
\newblock {\em Uspehi Mat. Nauk}, 32(4 (196)):55--112, 287, 1977.

\bibitem{PS}
Y.~Pesin and Y.~Sinai.
\newblock Gibbs measures for partially hyperbolic attractors.
\newblock {\em Ergodic Theory Dynam. Systems}, 2(3-4):417--438, 1982.

\bibitem{PotJMD}
R.~Potrie.
\newblock Partial hyperbolicity and foliations in $\mathbb{T}^3$.
\newblock {\em J. Mod. Dyn.} 9(8): 81--121, 2015.


\bibitem{Pot2015}
R.~Potrie.
\newblock A few remarks on partially hyperbolic diffeomorphisms of
  {$\Bbb{T}^3$} isotopic to {A}nosov.
\newblock {\em J. Dynam. Differential Equations}, 26(3):805--815, 2014.

\bibitem{PuSh89}
C.~Pugh and M.~Shub.
\newblock Ergodic attractors.
\newblock {\em Trans. Amer. Math. Soc.}, 312(1):1--54, 1989.

\bibitem{Hertz:2011vu}
F.~Rodriguez~Hertz, J.~Rodriguez~Hertz, A.~Tahzibi, and R.~Ures.
\newblock Uniqueness of {SRB} measures for transitive diffeomorphisms on
  surfaces.
\newblock {\em Comm. Math. Phys.}, 306(1):35--49, 2011.

\bibitem{HHU}
F.~Rodriguez~Hertz, J.~Rodriguez~Hertz, and R.~Ures.
\newblock A non-dynamically coherent example on $\mathbb{T}^3$.
\newblock {\em Ann. Inst. H. Poincar\'e Anal. Non Lin\'eaire}.
\newblock To appear.

\bibitem{U}
R.~Ures.
\newblock Intrinsic ergodicity of partially hyperbolic diffeomorphisms with a
  hyperbolic linear part.
\newblock {\em Proc. Amer. Math. Soc.}, 140(6):1973--1985, 2012.



\bibitem{W82}
Peter Walters.
\newblock {\em {An introduction to ergodic theory}}, volume~79 of {\em Graduate
  Texts in Mathematics}.
\newblock Springer-Verlag, New York, 1982.

\bibitem{Young:2002vc}
L.~Young.
\newblock {What are SRB measures, and which dynamical systems have them?}
\newblock {\em Journal of Statistical Physics}, 108(5-6):733--754, 2002.

\end{thebibliography}

\def\cprime{$'$}

\end{document}